\documentclass[11pt]{article}
\usepackage[utf8]{inputenc}
\usepackage{amsmath}
\usepackage{amssymb}
\usepackage{amsfonts}

\newtheorem{theorem}{Theorem}[section]

\newtheorem{example}[theorem]{Example}

\newtheorem{lemma}[theorem]{Lemma}

\newtheorem{remark}[theorem]{Remark}

\newenvironment{proof}[1][Proof]{\noindent\textbf{#1.} }{\ \rule{0.5em}{0.5em}}

\usepackage[numbers]{natbib}
\title{\centering On statistical estimation and inferences in optional regression models}
\author{ {Mohamed Abdelghani\thanks{Morgan Stanley, NY, USA.}, Alexander Melnikov\thanks{Department of Mathematical and Statistical Sciences, University of Alberta, Edmonton, Canada} and Andrey Pak\thanks{Department of Mathematical and Statistical Sciences, University of Alberta, Edmonton, Canada.
Corresponding Author: pak@ualberta.ca.}}}
\date{}
\date{March 2020}

\begin{document}

\maketitle

\begin{abstract}
    The main object of investigation in this paper is a very general regression model in optional setting – when an observed process is an optional semimartingale depending on an unknown parameter. It is well-known that statistical data may present an information flow/filtration without "usual conditions". The estimation problem is achieved by means of structural least squares (LS) estimates and their sequential versions. The main results of the paper are devoted to the strong consistency of such LS-estimates. For sequential LS-estimates the property of fixed accuracy is proved. 
\end{abstract}

\bigskip

\section{Introduction}
Regression Analysis is an integral part of Mathematical Statistics. Developments in this area are important from both theoretical and applied points of view. In statistics of random processes a regression model is considered as a semimartingale where the drift depends on an unknown parameter and the martingale part presents the errors in observations. Such a view point is very productive because it creates a possibility to study a variety of regression models (with discrete and continuous time) in an unified way, using martingale methods (see, for example, \cite{liptser2012theory}, \cite{mel1996stochastic}). 

The standard martingale theory is well-developed under so-called "usual conditions", when filtration (information flow) is complete and right-continuous. However, statistical data is usually delivered by a stochastic process, whose history (natural filtration) may not be right-continuous, and therefore such technical conditions may not be fulfilled (see \cite{abdmel2020}). This is the main reason why we need to consider regression models in more general setting which we call here the optional regression model. Optional semimartingales, on which our optional regression model is based, admit trajectories which are not right-continuous and arise when "usual conditions" are not assumed on filtered probability space. Up to our knowledge, currently there are no works devoted to the relaxing of these "usual conditions" and investigation of such general optional regression model. Whenever we use definitions or results from the theory of optional processes throughout the paper, we refer the reader to the book \cite{abdmel2020}.

In the first part of the paper, we focus on strong consistency of the proposed LS-estimate. In case of the observed process being cadlag (right-continuous with left limits) semimartingale this problem was extensively studied in \cite{mel1996stochastic} (also see \cite{mel1987regression}). 

In the second part of the paper, we concentrate our attention on the sequential estimates with guaranteed accuracy. In comparison to LS-estimates the sequential LS-estimates posses an advantage of having bounded variance.  This type of estimates in cadlag case is well-established (see  \cite{melnov88}, \cite{mel1996stochastic}, \cite{galtchouk2001sequential}).

The paper is organized in the following way: in section 2 we introduce the general regression model along with structural LS-estimates and auxiliary results. In section 3 we prove strong consistency of the proposed LS-estimates. In section 4 we consider sequential estimate, show that these estimates are unbiased and have a property of guaranteed accuracy under suitable conditions on regressor and error term. In addition, we investigate a problem related to hypothesis testing. Finally, in section 5 we present an extension of sequential LS-estimators for non-linear regression models and several illustrative examples.

\section{Stochastic regression model in optional setting}
Suppose that on the fixed stochastic basis $(\Omega, \mathcal{F}, \mathbf{F}=(\mathcal{F}_t)_{t\geq0}, \mathbf{P})$ without "usual conditions", we observe a one-dimensional process $X$. Let $\mathcal{O}$ and $\mathcal{P}$ be optional and predictable $\sigma $-algebras, respectively, as well as families of corresponding optional and predictable processes. 
In general, optional and predictable
processes have right and left limits but may not necessarily be right or left
continuous in $\mathbf{F}$. For either optional or predictable processes we
can define the following processes: $X_{-}=(X_{t-})_{t\geq 0}$ and $%
X_{+}=(X_{t+})_{t\geq 0}$, $\Delta X=(\Delta X_{t})_{t\geq
0} $ such that $\Delta X_{t}=X_{t}-X_{t-}$ and $\Delta
^{+}X=(\Delta ^{+}X_{t})_{t\geq 0}$ such that $\Delta
^{+}X_{t}=X_{t+}-X_{t}$.

We denote $\mathcal{P}_s$ as a collection of strongly predictable processes $a_t\in \mathcal{P}_s$, if $a_t\in\mathcal{P}$ and $a_{t+}\in \mathcal{O}.$  The families of increasing and increasing locally integrable processes are denoted by $\mathcal{V}^+$ and $\mathcal{A}^+_{loc}$, respectively. Let $\mathcal{M}$ and $\mathcal{M}^2\ (\mathcal{M}_{loc}$ and  $\mathcal{M}^2_{loc})$ denote the set of optional (local) martingales and optional square integrable (local) martingales, respectively.

Suppose the process $X$ has the following form
\begin{equation} \label{model}
    X_t=f\circ a_t\theta+M_t,
\end{equation}
where $f\circ a_t$ is an optional stochastic integral such that
$$f\circ a_t=\int_{]0,t]}f^r_sda^r_s +\int_{[0,t[}f^g_sda^g_{s+},$$
$a=a^r+a^g \in \mathcal{A}^+_{loc}\cap \mathcal{P}_s,\ M \in \mathcal{M}_{loc},\ f_t$ is a bilinear pair $f_t=(f^r_t,f^g_t),\ f^r_t \in \mathcal{P},\ f^g_t \in \mathcal{O},$
and $\theta\in \mathbb{R}$ is the unknown parameter which we need to estimate.

As the estimator of $\theta$ we consider the statistic
\begin{equation} \label{LSestimator}
    \theta_t=F^{-1}_t\left(f \circ X_t\right),
\end{equation}
where $F_t:=f^2 \circ a_t \in \mathcal{A}^+_{loc}\cap \mathcal{P}_s$
is assumed to be non-zero (a.s.). This assumption is not restrictive because further we suppose that $F_t \to \infty$ (a.s.) to provide strong consistency of $\theta_t.$

The structure of the estimator \eqref{LSestimator}
 is similar to estimator obtained by the method of Least Squares (LS) in classical regression analysis. Therefore, $\theta_t$ will be called the structural LS-estimator of $\theta$. It is well known how to study its asymptotic behaviour with the help of the Strong Law of Large Numbers (SLLN). Liptser (1980) \cite{liptser1980strong} proposed a very general form of SLLN for local martingales using a stochastic Kronecker's Lemma. For reader's convenience, let us reproduce this scheme in optional setting (see \cite{gasp}, \cite{mel1996stochastic}).

To prove Kronecker's Lemma in optional setting, we need the following result on sets of convergence of optional martingales.
In what follows, we denote $\widetilde{D}$ the compensator of some increasing process $D$.
\begin{lemma} \label{convset} 
If $Y \in \mathcal{M}_{loc}$ then 
$$(\widetilde{D}_\infty<\infty) \subseteq(Y\to)\ a.s.,$$ where 
\begin{equation*}
    D_t=\langle Y^c \rangle_t+\sum_{0<s\leq t}\frac{(\Delta Y_s)^2}{1+|\Delta Y_s|}+\sum_{0\leq s<t}\frac{(\Delta^+ Y_s)^2}{1+|\Delta^+ Y_s|},
\end{equation*}
and $(Y\to)$ is the set, on which there exists a finite random variable $Y_{\infty}(\omega)=\lim_{t\to \infty} Y_t(\omega)<\infty.$
\end{lemma}

Now we present the following generalization of Kronecker's Lemma.

\begin{lemma} \label{lemmaSLLN}
 For processes $N\in \mathcal{M}_{loc}$ and $A\in \mathcal{V}^+ \cap \mathcal{P}_s$ the following relation holds
$$(A_\infty=\infty)\cap (Y_t\rightarrow)\subseteq(A_t^{-1}N_t\rightarrow0) \ \  (a.s.)\ \  (t\rightarrow\infty),
$$
where 
\begin{equation}\label{ProcY}
Y_t=\int_{]0,t]}(1+A_s)^{-1}dN^r_s+\int_{[0,t[}(1+A_{s+})^{-1}dN^g_{s+}.
\end{equation}
\end{lemma}

\begin{proof}
 From \eqref{ProcY} it is easy to see that
\begin{equation}\label{Proc2}
    \int_{]0,t]}(1+A_s)dY^r_s+\int_{[0,t[}(1+A_{s+})dY^g_{s+}=N_t-N_0.
\end{equation}

Using integration by parts formula (see Lemma 3.4, \cite{amprisk2020}) we obtain
\begin{equation*}
    (1+A_t)Y_t=\int_{]0,t]}(1+A_s)dY^r_s+\int_{[0,t[}(1+A_{s+})dY^g_{s+}+\int_{]0,t]}Y_{s-}dA^r_s+\int_{[0,t[}Y_sdA^g_{s+}.
\end{equation*}
Then from \eqref{Proc2} we conclude that
\begin{equation*}
    \frac{N_t}{1+A_t}=\frac{N_0+Y_t}{1+A_t}+\frac{1}{1+A_t}(A_tY_t-\int_{]0,t]}Y_{s-}dA^r_s-\int_{[0,t[}Y_sdA^g_{s+}).
\end{equation*}

Since $\sup_{t\ge0} |Y_t| < \infty$ on the set $(A_\infty=\infty)\cap (Y\rightarrow)$, we have that $(1+A_t)^{-1}(N_0 +Y_t)\rightarrow 0$ a.s. as $t\rightarrow\infty$ .
On the other hand, we have for $u<t, v<t$
\begin{align} \notag
    \frac{1}{1+A_t}&\left|A_tY_t-\int_{]0,t]}Y_{s-}dA^r_s-\int_{[0,t[}Y_sdA^g_{s+}\right|\\ \notag
    =&\frac{1}{1+A_t}\left|\int_{]0,t]}(Y_t-Y_{s-})dA^r_s+\int_{[0,t[}(Y_t-Y_s)dA^g_{s+}\right|\\ \notag
     \leq& \frac{1}{1+A_t}\Bigg(\int_{]0,u]}|Y_t-Y_{s-}|dA^r_s+\int_{]u,t]}(|Y_\infty-Y_t|+|Y_\infty-Y_{s-}|)dA^r_s\\ \notag
    & +\int_{[0,v[}|Y_t-Y_s|dA^g_{s+}+\int_{[v,t[}(|Y_\infty-Y_t|+|Y_\infty-Y_{s}|)dA^g_{s+}\Bigg)\\ \notag
    \leq& 2\sup_{s\geq 0}|Y_s|(1+A_t)^{-1}(A^r_u+A^g_v)+|Y_\infty-Y_t|\\
    & +(1+A_t)^{-1}\int_{]u,t]}|Y_\infty-Y_{s-}|dA^r_s+(1+A_t)^{-1}\int_{[v,t[}|Y_\infty-Y_s|dA^g_{s+}.\label{big3}
\end{align}

Using the fact that $\sup_{s\geq0}|Y_s| <\infty$ on the set $(Y\rightarrow)$, we can choose for sufficiently large $t$ appropriate values $u$ and $v$ on the set $(A_\infty=\infty)\cap(Y\rightarrow)$ to make the right side of \eqref{big3} tend to zero.
Consequently, the statement of the lemma follows.
\end{proof}
\begin{remark}
 Although we proved Lemma 2.2 for the process $N \in \mathcal{M}_{loc}$, this proof also works for any optional semimartingale $N$.
\end{remark}
 
 \section{Strong Consistency}
In this section we will show that the estimator $\theta_t$ in \eqref{LSestimator} is strongly consistent. The proof of the strong consistency is based on the SLLN in optional case. 

Let $N\in \mathcal{M}_{loc}$ and
\begin{equation}\label{decomp}
N_t= N^c_t+\int_{]0,\infty]}\int_{\mathbb{R}_0}x d(\mu^r-\nu^r)_s+ \int_{[0,\infty[}\int_{\mathbb{R}_0}x d(\mu^g-\nu^g)_{s+}
\end{equation}
be the canonical decomposition of $N$, where $\mathbb{R}_0=\mathbb{R}\setminus \{0\},$ $N^c$ be a continuous part of $X$, $\mu^r$ and $\mu^g$ be random measures of right and left jumps of $N$, and $\nu^r$ and $\nu^g$ be their respective compensators. 
\begin{theorem} \label{SLLN}
Let $A\in \mathcal{V}^+ \cap\mathcal{P}_s$ and $A_\infty=\infty$ a.s. If $N\in \mathcal{M}_{loc}$ and for some $q\in [1,2]$
\begin{equation} \label{condSLLN}
\int_{]0,\infty]} \frac{d\langle N^c\rangle_s}{(1+A_s)^2}+ \int_{]0,\infty]}\int_{\mathbb{R}_0}|1+A_s|^{-q}|x|^qd \nu^r_s+ \int_{[0,\infty[}\int_{\mathbb{R}_0}|1+A_{s+}|^{-q}|x|^qd \nu^g_{s+} < \infty
\end{equation}
then
$$A_t^{-1}N_t\rightarrow0 \text{ a.s. as } t\rightarrow\infty.
$$
\end{theorem}

\begin{proof}
Using the fact that 
$$\frac{x^2}{1+|x|}\leq |x|^q, \  q\in [1,2],$$
we get for $q \in [1,2]$
\begin{align*}
    \Tilde{D}_\infty=&\langle Y^c \rangle_\infty + \widetilde{\sum_{s\leq \infty}\frac{(\Delta Y_s)^2}{1+|\Delta Y_s|}}+\widetilde{\sum_{s<\infty}\frac{(\Delta^+ Y_s)^2}{1+|\Delta^+ Y_s|}}\\
    \leq & \langle Y^c \rangle_\infty + \widetilde{\sum_{s\leq \infty}(\Delta Y_s)^q}+\widetilde{\sum_{s<\infty}(\Delta^+ Y_s)^q},
\end{align*}
where $Y$ as defined in \eqref{ProcY}.
Thus, from \eqref{condSLLN} it follows that
\begin{multline*}
    \Tilde{D}_\infty \leq \int_{]0,\infty]} \frac{d\langle N^c\rangle_s}{(1+A_s)^2}+ \int_{]0,\infty]}\int_{\mathbb{R}_0}|1+A_s|^{-q}|x|^qd \nu^r_s\\
    + \int_{[0,\infty[}\int_{\mathbb{R}_0}|1+A_{s+}|^{-q}|x|^qd \nu^g_{s+} < \infty.
\end{multline*}
By Theorem \ref{convset} and Lemma \ref{lemmaSLLN}, $ A_t^{-1}N_t\rightarrow0 \text{ a.s. as } t\rightarrow\infty$.
\end{proof}

\begin{theorem}
Suppose for the model \eqref{model} that $F_\infty=\infty$ and for some $q\in [1,2]$
\begin{multline}  
\int_{]0,\infty]} \frac {(f^r_s)^2d\langle M^c\rangle_s}{(1+F_s)^2}+ \int_{]0,\infty]}\int_{\mathbb{R}_0}|1+F^r_s|^{-q}|f^r_s|^q|x|^qd \nu^r_s\\
+ \int_{[0,\infty[}\int_{\mathbb{R}_0}|1+F^g_{s}|^{-q}|f^g_s|^q|x|^qd \nu^g_{s+} < \infty.  \label{condconsis}
\end{multline}
Then $\theta_t \to \theta$ (a.s.) as $t \to \infty.$
\end{theorem}
\begin{proof}
It is sufficient to note that 
$$\theta_t-\theta=A^{-1}_tN_t,$$
where $A_t:=F_t$ and $N_t:=f\circ M_t$.
 By Theorem \ref{SLLN} we get immediately the statement of the theorem.
\end{proof}

\section{Sequential LS-estimators}
Let us consider the model \eqref{model} with $M\in \mathcal{M}^2_{loc} (\mathbb{R})$. We assume that there exists a non-negative random variable $\xi$ such that

\begin{equation}\label{assumption1}
    \frac{d{\langle M \rangle_t}}{da_t} \leq \xi , \ \ \ \ f^2 \circ a_t \in \mathcal{A}^+_{loc}\cap \mathcal{P}_s.
\end{equation}

Next, for fixed $H$ we define
\begin{equation} \label{tauH}
    \tau_H =\inf\{t:f^2 \circ a_t\geq H\},
\end{equation}
with $\tau_H=\infty$ if the corresponding set is empty. 
We assumed that processes $f^r\in \mathcal{P},\ f^g \in \mathcal{O}$ and $a \in \mathcal{P}_s$, consequently, $\{t: f^2 \circ a_t\geq H\}$ are stopping times. Therefore, by Theorem 2.4.14 in \cite{abdmel2020}  $\tau_H$ is a wide sense stopping time.

On the set $\{\tau_H < \infty\}$ we define a random variable $\beta_H$ by the relation 
\begin{equation} \label{H}
    f^2 \circ a_{\tau_H-} + \beta_H \left( (f^r_{\tau_H})^2\Delta a_{\tau_H}+ (f^g_{\tau_H})^2\Delta^+ a_{\tau_H}\right)=H,
\end{equation} 
and on $\tau_H= \infty$ we put $\beta_H=0.$ Then $\beta_H \in [0,1]$ and is a $\mathcal{F}_{\tau_H}$-measurable random variable.

We consider the following statistic as an estimator of $\theta$
\begin{equation}\label{estimator}
    \hat{\theta}_H=H^{-1}\left[  f^2 \circ X_{\tau_H-} + \beta_H \left( (f^r_{\tau_H})^2\Delta X_{\tau_H}+ (f^g_{\tau_H})^2\Delta^+ X_{\tau_H}\right) \right].
\end{equation}

The next theorem shows that the statistic defined by means of \eqref{assumption1}-\eqref{estimator} is an unbiased estimator of $\theta$ and has the property of guaranteed accuracy, i.e., bounded variance.

\begin{theorem}\label{seqest}
Suppose that assumptions \eqref{assumption1} hold, $\mathbf{E}\xi < \infty,$ and
\begin{equation}\label{assumption2}
    \mathbf{P}\{ f^2 \circ a_\infty=\infty\}=1.
\end{equation}
Then for all $H>0$
$$\mathbf{P}\{\tau_H<\infty\}=1, \ \ \ \mathbf{E} \hat{\theta}_H=\theta, \ \ \ \ \mathbf{Var}\hat{\theta}_H\leq H^{-1}\mathbf{E} \xi.$$
\end{theorem}

\begin{proof}
First,
\begin{align*}
    \mathbf{P}\{\tau_H=\infty\}&=\mathbf{P}\{ f^2 \circ a_\infty <H\}\\
    &=1-\mathbf{P}\{ f^2 \circ a_\infty \geq H\}\\
    &\leq 1-\mathbf{P}\{ f^2 \circ a_\infty =\infty\}\\
    &=1-1=0.
\end{align*}
Thus, 
$$\mathbf{P}\{\tau_H<\infty\}=1-\mathbf{P}\{\tau_H=\infty\}=1.$$
Next, using \eqref{model} and \eqref{estimator}, we obtain
\begin{align*}
        \hat{\theta}_H=&H^{-1}\left[  f^2 \circ X_{\tau_H-} + \beta_H \left( (f^r_{\tau_H})^2\Delta X_{\tau_H}+ (f^g_{\tau_H})^2\Delta^+ X_{\tau_H}\right) \right]\\
        =&H^{-1}\big[  (f^2 \circ a_{\tau_H-})\theta + \beta_H \left( (f^r_{\tau_H})^2\Delta a_{\tau_H}+ (f^g_{\tau_H})^2\Delta^+ a_{\tau_H}\right)\theta \\
        &+ f^2 \circ M_{\tau_H-} + \beta_H \left( (f^h_{\tau_H})^2\Delta M_{\tau_H}+ (f^g_{\tau_H})^2\Delta^+ M_{\tau_H}\right)\big]\\
        =&\theta+H^{-1}N_{\tau_H},
\end{align*}
where 
$$N_t=I_{\{t<\tau_H\}} f^2 \circ M_t + I_{\{t=\tau_H\}}\beta_H \left((f^r_{\tau_H})^2\Delta M_{\tau_H}+ (f^g_{\tau_H})^2\Delta^+ M_{\tau_H}\right)$$

Since the process $N_t$ is a stochastic integral with respect to the optional square integrable local martingale $M$, by the properties of optional stochastic integrals we have 
$$\langle N \rangle_t=I_{\{t<\tau_H\}}f^2 \circ \langle M \rangle_t + I_{\{t=\tau_H\}}\beta^2_H \left( (f^r_{\tau_H})^2\Delta \langle M \rangle_{\tau_H}+ (f^g_{\tau_H})^2\Delta^+ \langle M \rangle_{\tau_H}\right).$$
Hence by \eqref{assumption1} and \eqref{H} we get
$$\langle N \rangle_{\tau_H}\leq \xi \left[ f^2 \circ a_{\tau_H-} + \beta_H \left( (f^r_{\tau_H})^2\Delta a_{\tau_H}+(f^g_{\tau_H})^2\Delta^+ a_{\tau_H}\right)\right]=\xi H.$$
Consequently, $N_{t\wedge \tau_h}$ is an optional square integrable martingale, and therefore 
$$\mathbf{E}N_{\tau_H}=0, \ \ \ \mathbf{E}N^2_{\tau_H}\leq H\mathbf{E} \xi,$$
which proves the theorem.
\end{proof}

Now, it is reasonable to discuss the following problem of distinguishing two hypotheses with simultaneous estimation of the parameter $\theta \in \mathbb{R}:$

\begin{align}\label{hypo}
    H_0:\ X_t=& f\circ a_t \theta + M_t,\\
    H_1:\ X_t=&M_t,
\end{align}
where $M\in \mathcal{M}^2_{loc}(\mathbb{R}).$
Assuming that \eqref{assumption1} is fulfilled, we define $\tau_H, \ \beta_H$ as in \eqref{tauH}, \eqref{H} and 
$$\phi_H(X)=H^{-1}\left[  f^2 \circ X_{\tau_H-} + \beta_H \left( (f^r_{\tau_H})^2\Delta X_{\tau_H}+ (f^g_{\tau_H})^2\Delta^+ X_{\tau_H}\right) \right].$$
\begin{theorem}
Suppose that in the problem \eqref{hypo} the parameter satisfies 
$$\theta\in \ R_\delta=\{\theta\in \mathbb{R}: |\theta|\geq \delta >0\}$$
and condition \eqref{assumption2} is fulfilled for both hypotheses $H_0$ and $H_1, \ \mathbf{E} \xi <\infty.$ Then for a given $\epsilon>0$ the criterion 

$$\Delta (\tau_H)= 
\begin{cases}
 0 \text{ if } |\phi_H(X)|\geq \delta/2,\\    
  1 \text{ if } |\phi_H(X)|< \delta/2,
\end{cases}$$
for $H\geq 4 \delta^{-2}\epsilon^{-1}\mathbf{E}\xi$ ensures distinguishability of the hypotheses $H_0$ and $H_1$ with probabilities of errors not exceeding $\epsilon>0.$
\end{theorem}
\begin{proof}
From the definition of $\phi_H(X)$, if either $H_0$ or $H_1$ is true, then it follows, respectively, that
$$\phi_H(X)=\theta+ \phi_H(M) \ \  \text{  or  } \ \  \phi_H(X)=\phi_H(M).$$
Note that by Theorem 4.1 under either of the hypotheses $H_0,\  H_1$
\begin{equation}\label{estimates}
    \mathbf{E}\phi_H(M)=0, \ \ \mathbf{E}\phi^2_H(M)\leq H^{-1}\mathbf{E}\xi.
\end{equation}
Then in case when $H_0$ is true, applying Chebyshev's inequality we obtain that for $4\delta^{-2}\epsilon^{-1}\mathbf{E}\xi \leq H$
\begin{align*}
    \mathbf{P}\{\omega: \Delta(\tau_H)\ne 1\}=&\mathbf{P}\{\omega: |\phi_H(X)|\geq \delta/2\}\\
    =&\mathbf{P}\{\omega:|\phi_H(M)|\geq \delta/2\|\leq 4\delta^{-2}\mathbf{E}|\phi_H(M)|^2\\
    \leq&4\delta^{-2}H^{-1}\mathbf{E}\xi\leq \epsilon.
\end{align*}
In case when $H_1$ is true, using the simple fact that for $\theta \in R_\delta$
$$\delta-|\phi_H(M)|\leq|\delta+\phi_H(M)|\leq |\theta+\phi_H(M)|=|\phi_H(X)|$$
implying that
 $$\{\omega:|\phi_H(X)|<\delta/2\} \subseteq \{\omega: |\phi_H(M)|\geq \delta/2\},$$
 we arrive to the following estimate of the probability of error:
 \begin{align*}
     \sup_{\theta\in R_\delta} \mathbf{P}\{\omega: \Delta (\tau_H)\ne 0\}=&\sup_{\theta\in R_\delta}\mathbf{P}\{\omega: |\phi_H(X)|<\delta/2\}\\
     \leq& \sup_{\theta\in R_\delta}\mathbf{P}\{\omega: |\phi_H(M)|<\delta/2\}\\
     \leq& 4\delta^{-2}\sup_\theta \mathbf{E}\phi^2_H(M)\leq 4\delta^{-2}H^{-1}\mathbf{E}\xi\leq \epsilon.
 \end{align*}
\end{proof}

\bigskip
\section{Further extensions and examples}
Let us show how the linear optional regression model \eqref{model} can be extended in a non-linear case (see, for example, \cite{bishwal2006}). The non-linear optional regression model has the following form 
\begin{equation}
    X_t=f\circ a_t g(\theta)+ M_t,
\end{equation}
where the processes $f, a, M$ satisfy the same assumptions as in \eqref{assumption1}-\eqref{H} and $g: \mathbb{R}\to \mathbb{R}$ is a continuous, bijective function with a continuous inverse $g^{-1}$. Let $\tau_H$ and $\beta_H$ be as in \eqref{tauH} and \eqref{H}, respectively. Then the sequential LS-estimator can be obtained by defining $\zeta:=g(\theta)$ and realizing from \eqref{estimator} that
\begin{align}\label{nonlinestimator}\notag
    \hat{\zeta}_H&=g(\hat{\theta}_H)=H^{-1}\left[  f^2 \circ X_{\tau_H-} + \beta_H \left( (f^r_{\tau_H})^2\Delta X_{\tau_H}+ (f^g_{\tau_H})^2\Delta^+ X_{\tau_H}\right) \right] \text{  or }\\
    \hat{\theta}_H&=g^{-1}\left(H^{-1}\left[  f^2 \circ X_{\tau_H-} + \beta_H \left( (f^r_{\tau_H})^2\Delta X_{\tau_H}+ (f^g_{\tau_H})^2\Delta^+ X_{\tau_H}\right) \right]\right).
\end{align}

Now, suppose \eqref{assumption2} holds, $g(\theta)$ is differentiable and 
\begin{equation}\label{gcond}
    \int^\infty_{-\infty}(g^{-1}(x))^2\exp(-x^2/2)dx<\infty.
\end{equation}
Using the same argument as in the proof of Theorem \ref{seqest}, we show that $\mathbf{P}\{\tau_H=\infty\}=1.$

Next, note that 
\begin{align*}
    \mathbf{E}(\theta_H-\theta)=\mathbf{E}\left[g^{-1}(g(\theta)+H^{-1}N_{\tau_H})-g^{-1}(g(\theta))\right],\\
    \mathbf{E}(\theta_H-\theta)^2=\mathbf{E}\left[g^{-1}(g(\theta)+H^{-1}N_{\tau_H})-g^{-1}(g(\theta))\right]^2,
\end{align*}
where 
$$N_t=I_{\{t<\tau_H\}} f^2 \circ M_t + I_{\{t=\tau_H\}}\beta_H \left((f^r_{\tau_H})^2\Delta M_{\tau_H}+ (f^r_{\tau_H})^2\Delta^+ M_{\tau_H}\right).$$
From the proof of Theorem \ref{seqest} we already know that $N_{t\wedge \tau}\in \mathcal{M}\cap \mathcal{M}^2, \ \mathbf{E}N_{\tau_H}=0$ and $\mathbf{E}N^2_{\tau_H}\leq H \mathbf{E}\xi.$ From assumption \eqref{assumption2} it follows that $\langle N\rangle_\infty=\infty.$
Thus, by Theorem \eqref{SLLN} we have 
\begin{equation*}
    \lim_{t\to \infty}\frac{N_t}{\langle N \rangle_t}=0 \  (a.s.) \text{  and } \lim_{H\to \infty}\frac{N_{\tau_H}}{\langle N \rangle_{\tau_H}}=0 \ (a.s.).
\end{equation*}

Since $\langle N \rangle_{\tau_H}\leq H\xi$, we get
$$\lim_{H\to \infty}\frac{N_{\tau_H}}{H}=0 \ (a.s.).$$
Using the Skorokhod embedding theorem we obtain
a\begin{align*}
    \mathbf{E}(\theta_H-\theta)=&\int^{\infty}_{-\infty}\left[g^{-1}(g(\theta)+H^{-1}x)-g^{-1}(g(\theta))\right]e^{-x^2/(2H)}dx\\
    =&\int_{y\leq A}\left[g^{-1}(g(\theta)+H^{-1/2}y)-g^{-1}(g(\theta))\right]e^{-y^2/2}dy\\
    &+\int_{y> A}\left[g^{-1}(g(\theta)+H^{-1/2}y)-g^{-1}(g(\theta))\right]e^{-y^2/2}dy
\end{align*}
Applying condition \eqref{assumption2} and \eqref{gcond} one can always choose numbers $A_0(\epsilon)$ and $H_0(\epsilon)$ such that for $H\geq H_0(\epsilon)$
$$\int_{y\leq A_0(\epsilon)}\left[g^{-1}(g(\theta)+H^{-1/2}y)-g^{-1}(g(\theta))\right]e^{-y^2/2}dy< \epsilon$$
and 
$$\int_{y> A_0(\epsilon)}\left[g^{-1}(g(\theta)+H^{-1/2}y)-g^{-1}(g(\theta))\right]e^{-y^2/2}dy< \epsilon.$$
Hence, $\mathbf{E}(\theta_H-\theta)< 2 \epsilon$ for $H> H_0(\epsilon).$
This proves $\theta_H$ is asymptotically unbiased as $H\to \infty.$

It is not difficult to show similar calculations for $\mathbf{E}(\theta_{\tau_H}-\theta)^2.$ From Cramer-Rao-Wolfovitz inequality it follows that $\theta_{\tau_H}$ is asymptotically efficient.

On the other hand, we can show that $\theta_H$ is unbiased and efficient estimator under assumption that $g^{-1}(\theta)$ is differentiable and has a bounded first derivative. Using Mean Value theorem, we have
\begin{align*}
      \mathbf{E}|\theta_H-\theta|=&\mathbf{E}\left[g^{-1}(g(\theta)+H^{-1}N_{\tau_H})-g^{-1}(g(\theta))\right]\\
      \leq &\mathbf{E}\sup_{\zeta\in[g(\theta),g(\theta)+H^{-1}N_{\tau_H}]}|(g^{-1})' (\zeta)|\frac{N_{\tau_H}}{H}\\
      \leq&0
\end{align*}
Similarly,
\begin{align*}
      \mathbf{E}|\theta_H-\theta|^2 \leq &\mathbf{E}\left[\sup_{\zeta\in[g(\theta),g(\theta)+H^{-1}N_{\tau_H}]}|(g^{-1})' (\zeta)|\right]^2\frac{N^2_{\tau_H}}{H^2}\\
      \leq&K^2 H^{-1}\mathbf{E}\xi,
\end{align*}
where $K$ is the constant bound on $(g^{-1})'(\zeta).$

Let us now illustrate several examples.
\begin{example}\textbf{Non-linear regression model.}
 Consider the following non-linear model
 $$X_t=f\circ a_t\sqrt{\theta}+ M_t,$$
 where $f, a, M $ satisfy assumption \eqref{assumption1}-\eqref{H},\eqref{assumption2}. 
 
 The function $g(\theta)=\sqrt{\theta}, \theta\geq 0,$ is differentiable and its inverse $g^{-1}(\theta)=\theta^2$ clearly satisfies \eqref{gcond}:
 $$    \int^\infty_{-\infty}(g^{-1}(x))^2\exp(-x^2/2)dx=    \int^\infty_{-\infty}x^4\exp(-x^2/2)dx=3\sqrt{2\pi}<\infty.$$
 Thus, by the above discussion the sequential LS-estimator $$\hat{\theta}_H=\left(H^{-1}\left[  f^2 \circ X_{\tau_H-} + \beta_H \left( (f^r_{\tau_H})^2\Delta X_{\tau_H}+ (f^g_{\tau_H})^2\Delta^+ X_{\tau_H}\right) \right]\right)^2$$
 is both asymptotically unbiased and asymptotically efficient as $H\to \infty.$ 
 
\end{example}

\begin{example}\textbf{Risk process.}
Consider the following risk process 
\begin{equation}
    X_t=ct+\sigma W_t - a N^r_t +b N^g_t,
\end{equation}
where $c, \sigma, a, b$ are positive constants, $W$ is a Wiener process, $N^r$ and $N^g$ are a Poisson process and left-continuous modification of a Poisson process, respectively. The constant $c$ usually represents premium payments in risk theory, whereas $a$ and $b$ describe average value of claims and positive gains, respectively. The process $\sigma W_t$ is a random perturbation. 

We can rewrite the process $X_t$ as follows
\begin{equation*}
    X_t=\theta t+M_t,
\end{equation*}
where $\theta :=c-a\lambda^r+b\lambda^g, \  M_t:=\sigma W_t- a(N^r_t-\lambda^rt)+b(N^g_t-\lambda^gt),$ $\lambda^r$ and $\lambda^g$ are jump intensities of $N^r_t$ and $N^g_t$, respectively. 

The structural LS estimator of $\theta$ is
\begin{equation} \label{example1}
    \theta_t= \frac{X_t}{t}.
\end{equation}
The condition \eqref{condconsis} of Theorem 3.2, i.e.,
$$\sigma^2\int_{]0,\infty]} \frac{ds}{(1+s)^2}+ a^q\lambda^r\int_{]0,\infty]}(1+s)^{-q}ds< \infty,$$
holds for any $q\in (1,2]$. Thus, the estimator \eqref{example1} is strongly consistent.

The sequential LS-estimators have the following form
$$\hat{\theta}_H=\frac{X_{\tau_H}}{H}.$$
All assumptions of Theorem \ref{seqest} are obviously satisfied, and $\xi =\sigma^2+a^2\lambda^r+b^2\lambda^g$. Consequently,
$$\mathbf{P}\{\tau_H<\infty\}=1, \ \ \ \mathbf{E}\hat{\theta}_H=\theta, \ \ \ \mathbf{Var}\hat{\theta}_H\leq \frac{\sigma^2+a^2\lambda^r+b^2\lambda^g}{H}.$$
We usually want the process $X_t$ to be positive, so the estimator $\theta_t=c-a\lambda^r+b\lambda^g >0.$ This assertion is called as a net profit condition in risk theory.

\end{example}

\begin{example} \textbf{Ornstein–Uhlenbeck process.}
In mathematical finance we often deal with Ornstein–Uhlenbeck type processes that possess mean reversion property, i.e.

\begin{equation*}
    X_t= \int^t_0(\mu-X_{s-})ds\theta+M_t,
\end{equation*}
where $\mu $ is a positive constant, and $M \in \mathcal{M}^2_{loc}$.

We assume that 
\begin{equation}
    F_t:=\int^t_0(\mu-X_{s-})^2ds\in \mathcal{A}^+_{loc}\cap \mathcal{P}_s \ \text{ and } \ \frac{d\langle M\rangle_t}{dt}\leq \xi.
\end{equation}
Then, the structural LS estimator of $\theta$ is
\begin{equation} \label{example2}
    \theta_t= \frac{(\mu-X_-) \circ X_t}{\int^t_0(\mu-X_{s-})^2ds},
\end{equation}
and sequential LS-estimators have the following form
$$\hat{\theta}_H=H^{-1}\left[(\mu-X_{s-})^2\circ X_{\tau_H-}\right].$$

If the following condition
\begin{equation*}
   \int_{]0,\infty]} \frac{(\mu-X_{s-})^2d\langle M^c\rangle_s}{(1+F_s)^2}+ \int_{]0,\infty]}\int_{\mathbb{R}_0}|1+F^r_s|^{-q}|\mu-X_{s-}|^q|x|^qd \nu^r_s< \infty 
\end{equation*}  
holds for some $q\in [1,2]$ and $F_\infty=\infty$, then by Theorem 3.2 the estimator \eqref{example2} is strongly consistent. 

Furthermore, by Theorem \ref{seqest}, we have
$$\mathbf{P}\{\tau_H<\infty\}=1, \ \ \ \mathbf{E}\hat{\theta}_H=\theta, \ \ \ \mathbf{Var}\hat{\theta}_H\leq H^{-1}\mathbf{E}\xi.$$

\end{example}

\begin{example}
Finally, consider a regression model with well-known centered Gaussian martingale $M \in \mathcal{M}^2_{loc}$  and a deterministic function $f_t$,
\begin{equation}
    X_t= \int^t_0 f_sds\theta + M_t.
\end{equation}
Then, the LS estimator of $\theta$ is
\begin{equation} \label{example3}
    \theta_t= \frac{f \circ X_t}{\int^t_0f^2_sds},
\end{equation}
It can be shown in the same way as in \cite{mel1996stochastic} that strong consistency of \eqref{example3} follows only from the assumption of $\int^\infty_0f^2_sds=\infty$.

We assume that 
\begin{equation}
    \frac{d\langle M\rangle_t}{dt}\leq \xi,
\end{equation}
where $\xi$ is constant. Note that in case of centered Gaussian martingales $\langle M\rangle_t=\mathbf{E}M^2_t<\infty$ is a deterministic function.
Then sequential LS-estimators have the following form
$$    \hat{\theta}_H=H^{-1}\left[  f^2 \circ X_{\tau_H-} \right].$$
If, in addition, 
\begin{equation*}
    \mathbf{P}\left\{\int^\infty_0f^2_sds=\infty\right\}=1,
\end{equation*}
then, by Theorem \ref{seqest}, we have
$$\mathbf{P}\{\tau_H<\infty\}=1, \ \ \ \mathbf{E}\hat{\theta}_H=\theta, \ \ \ \mathbf{Var}\hat{\theta}_H\leq H^{-1}\xi.$$
\end{example}

\textbf{Funding}. The research was supported by the NSERC Discovery Grant RES0043487.

\bigskip
\textbf{Conflict of Interest Statement}. This paper is the private opinion of the authors and does not necessarily reflect the policy and views of Morgan Stanley.

\newpage

\bibliography{ref}

\end{document}